\documentclass[12pt]{article}
\usepackage{amsmath, amsfonts, amssymb, amsthm} 
\newtheorem{theorem}{Theorem}[section]
\newtheorem{proposition}[theorem]{Proposition}
\newtheorem{definition}[theorem]{Definition}
\newtheorem{remark}[theorem]{Remark}

\newtheorem{example}[theorem]{Example}

\newtheorem{corollary}[theorem]{Corollary}

\def\R{\mathbb R} \def\Z{\mathbb Z} \def\C{\mathbb C} 

\def\N{\mathbb N}
\def\C{{\mathbb C}} 
\def\Q{\mathbb Q}

\def\<{\,<\!}
\def\>{\!>\,}
\usepackage{graphics}
\usepackage{epsf}
\usepackage{epsfig}
\begin{document}

\title{Continued fraction expansions  for complex numbers - a general approach} 

\author{S.G. Dani}

\date{}

\maketitle

\renewcommand{\thefootnote}{}

\footnote{2010 \emph{Mathematics Subject Classification}: Primary 11J70; 
Secondary 11J25.}

\footnote{\emph{Key words and phrases}: Continued fraction expansions of 
complex numbers, algorithms, 
Eisenstein integers, Lagrange theorem, growth of denominators of convergents. }

\renewcommand{\thefootnote}{\arabic{footnote}}
\setcounter{footnote}{0}

\begin{abstract}

We introduce here a general framework for studying continued fraction 
expansions for complex numbers and establish 
some results on the convergence
of the corresponding sequence of convergents. For continued fraction 
expansions  with
partial quotients in a  discrete subring of $\C$
an analogue of the 
classical Lagrange theorem, characterising quadratic surds as numbers 
with eventually 
periodic continued fraction expansions, is proved. Monotonicity and exponential 
growth are established for the absolute
values of the denominators of the convergents for a class of continued
fraction algorithms with partial quotients in the ring of Eisenstein 
integers. 
\end{abstract}

\section{Introduction} 

A. Hurwitz \cite{H} introduced, in 1887, continued fraction expansions 
for complex numbers 
with Gaussian integers as partial quotients, via the nearest 
integer algorithm (known subsequently also as Hurwitz algorithm) 
and established some basic
properties concerning convergence of the sequence of convergents, 
and also proved an analogue of the 
classical Lagrange theorem characterizing quadratic surds as the
numbers with   
eventually periodic continued fractions; analogous results were also 
proved for the nearest integer algorithms with respect to Eisenstein integers
 as partial quotients, in place of Gaussian integers.

Application of complex continued fractions, typically involving the nearest
integer algorithm, to questions in Diophantine approximation
analogous to the theory for simple continued fractions for real numbers, 
was taken up by various authors during the last century (see \cite{Lv}, 
\cite{P}, \cite{L}, \cite{Hen},
and other references cited therein). 

In \cite{DN-compl}, where we considered  the  question of values of 
 binary quadratic forms with complex coefficients over pairs of Gaussian
 integers, we  extended the study  
of continued fractions to other possible algorithms in place of the nearest integer algorithm, 
and also introduced certain non-algorithmic constructions for continued fraction
expansions, via what were called iteration sequences; the partial quotients
for the continued fractions were however retained to be Gaussian integers. 
In this paper we set up a broader framework for studying continued 
fraction expansions for complex numbers, and prove certain 
general results on convergence, analogue of the Lagrange theorem, speed 
of convergence etc.. Our 
results   in particular generalize those of Hurwitz in the case of the nearest 
integer algorithms with respect to 
Gaussian integers and Eisenstein integers. 

\section{Preliminaries on continued fraction expansions}

We begin with a general formulation of the notion of continued fraction 
expansion, with flexible choices for the partial quotients. 
Let $\C$ denote the field of complex numbers and   $\C^*$ the set of nonzero 
numbers in $\C$. When $z\in \C$ can be expressed as 
$$z=a_0+
\displaystyle{\frac 1{a_1+\frac 1{a_2+\cdots}}},$$
with $a_j\in \C^*$ for all $j\in \N$ (natural numbers), where the 
right hand side is assigned the usual meaning as the limit of the truncated 
expressions (assuming that they represent genuine 
complex numbers and the limit exists - see below), we consider the 
expression  as above to be a  
continued fraction expansion for $z$; though our main application 
will be with $a_n$'s in specific rings, we shall first 
discuss some results
in which $a_n$ can be more general complex numbers. The above concept 
can be formulated more systematically as follows. 

 Let $\{a_n\}_{n=0}^\infty$ be a sequence
in $\C^*$. We associate to it two
sequences $\{p_n\}_{n=-1}^\infty$ and $\{q_n\}_{n=-1}^\infty$  defined 
recursively by the relations
$$ p_{-1}=1, p_0=a_0, p_{n+1}=a_{n+1}p_n+p_{n-1}, \; \mbox{for all }
n\geq 0, \mbox{ and}
$$
$$
 q_{-1}=0, q_0=1, q_{n+1}=a_{n+1}q_n+q_{n-1}, \; \mbox{
  for all } n \geq 0.
$$
If $q_n\neq 0$ for all $n$ then we can form
$p_n/q_n$, and if they converge, as $n\to \infty$, to a complex number
$z$ we say that {\it $\{a_n\}_{n=0}^\infty$ defines a continued
  fraction expansion of $z$}; in this case we express  $z$ as $[a_0,a_1, 
\dots]$. 

In conformity with the nomenclature adopted in $\cite{DN-compl}$ we 
call $\{p_n\}, \{q_n\}$ the {\it $\cal Q$-pair of sequences associated
to $\{a_n\}_{n=0}^\infty$}; ($\cal Q$ signifies ``quotient'').  The ratios $p_n/q_n$, with
$q_n\neq 0$ are called the {\it convergents} corresponding to the 
$\cal Q$-pair, or the sequence  $\{a_n\}_{n=0}^\infty$. 
We note that $p_nq_{n-1}-q_np_{n-1}=(-1)^{n-1}$ 
for all $n\geq 0$, as may be verified inductively.

\medskip

Given a $z\in \C^*$ ``candidates'' for continued fraction expansions for $z$ 
can be
arrived at by setting $a_n=z_n-z_{n+1}^{-1}$ for all $n\geq 0$, 
where $\{z_n\}_{n=0}^\infty$
is a sequence in $\C^*$ such that $z_0=z$ and for all $n\geq 1$, $|z_n|\geq 1$ 
and $z_{n+1}\neq z_n^{-1}$.
We shall call such a sequence an {\it iteration sequence for} $z$, and 
$\{a_n\}_{n=0}^\infty$ the associated sequence of {\it partial quotients}. (In 
\cite{DN-compl} ``iteration sequences'' were introduced, with slightly 
different conditions,
and $a_n$'s restricted to Gaussian integers). Whether a sequence of partial
quotients so constructed 
indeed defines a continued fraction expansion for $z$ is an issue that needs
to be considered however. 

We begin by noting the following general properties. 

\begin{proposition}\label{prop}
Let $z\in \C^*$ and let $\{z_n\}$ be an iteration sequence for $z$. Let 
$\{a_n\}_{n=0}^\infty$ be the associated sequence of partial
quotients, and let  $\{p_n\}, \{q_n\}$
be the $\cal Q$-pair of sequences associated to $\{a_n\}$. Then 
 for all $n\geq 0$ the following statements hold:

i)  $q_nz -p_n =(-1)^n(z_1\cdots z_{n+1})^{-1}$; 

ii) if $|p_n|>|z_1|^{-1}$ then $q_n\neq 0$;

iii) $(z_{n+1}q_n+q_{n-1})z=z_{n+1}p_n+p_{n-1}$; 

iv) if  $|q_{n-1}|<|q_n|$,
$|z-\frac{p_n}{q_n}|\leq |q_n|^{-2}(|z_{n+1}|-|\frac{q_{n-1}}{q_n}|)^{-1}$;  

v) if $q_n$'s are nonzero and $|q_n|\to \infty$ then 
${p_n}/{q_n}$ converges to $ z$ as $n\to \infty$. 

\end{proposition}

\proof i) We  argue by induction. Note that as $p_0=a_0$, $q_0=1$
and $z-a_0=z_1^{-1}$, the statement holds for $n=0$. Now let $n\geq 1$
and suppose by induction that the assertion holds for $0, 1, \dots, n-1$. 
Then we have 
$q_nz-p_n=(a_nq_{n-1}+q_{n-2})z-(a_np_{n-1}+p_{n-2})=a_n(q_{n-1}z-p_{n-1})+
(q_{n-2}z-p_{n-2})= (-1)^{n-1}(z_1\cdots 
z_n)^{-1}a_n+(-1)^{n-2}(z_1\cdots z_{n-1})^{-1}=(-1)^n(z_1\cdots 
z_n)^{-1}(-a_n+z_n)=(-1)^n(z_1\cdots z_{n+1})^{-1} $, which proves~(i).

ii) For $n\geq 0$ if $|p_n|>|z_1|^{-1}$, then by (i) we have $|q_nz|\geq |p_n|
-|z_1\cdots z_{n+1}|^{-1}\geq |p_n|-|z_1|^{-1}>0$, and hence $q_n\neq 0$, which 
proves~(ii). 

iii) For $n\geq 0$, by (i) we have $z_{n+1}(q_nz-p_n)=(-1)^n
(z_1\cdots z_{n+1})^{-1}z_{n+1}=(-1)^n(z_1\cdots z_{n})^{-1}=-(q_{n-1}z-p_{n-1})$, 
and hence $(z_{n+1}q_n+q_{n-1})z=z_{n+1}p_n+p_{n-1}$, which proves~(iii).

iv) By (iii) we get $|(z_{n+1}q_n+q_{n-1})(q_nz
-p_n)|=|(z_{n+1}p_n+p_{n-1})q_n-(z_{n+1}q_n+q_{n-1})p_n)| 
=|p_{n-1}q_n -q_{n-1}p_n|=1$. Also,  $|z_{n+1}q_n+q_{n-1}|\geq   
|q_n|(|z_{n+1}|-|\frac {q_{n-1}}{q_n}|)$, and we note that since $|z_{n+1}|\geq 1$
and $|q_{n-1}|<|q_n|$ we have $|z_{n+1}|-|\frac {q_{n-1}}{q_n}|>0$.
Thus $|z-\frac {p_n}{q_n}|=|q_n|^{-1}|z_{n+1}q_n+q_{n-1}|^{-1}\leq |q_n|^{-2}
(|z_{n+1}|-|\frac {q_{n-1}}{q_n})^{-1}$, which proves~(iv). 

v) If $q_n$ are nonzero and $|q_n|\to \infty$ then
$|z-\frac{p_n}{q_n}|=|q_n|^{-1}|z_1\cdots z_{n+1}|^{-1}\leq |q_n|^{-1}\to 0$, 
 and hence ${p_n}/{q_n}$ converges to $z$ as $n\to \infty$. This proves~(v). 
\hfill $\Box$

\medskip
We next specialise to sequences $\{a_n\}_{n=0}^\infty$ contained in 
discrete subrings of $\C$; by a subring we shall always mean one 
containing $1$, the multiplicative identity. When  $\{a_n\}_{n=0}^\infty$ is 
contained
in a discrete subring $\Gamma$, from the recurrence relations it 
follows that for the corresponding  $\cal Q$-pair $\{p_n\}, \{q_n\}$, 
we have $p_n,q_n\in  \Gamma$ for all~$n$.

\begin{proposition}\label{prop2}
Let the notation be as in Proposition~\ref{prop} and suppose further
that 

i) $\{a_n\}_{n=0}^\infty$ is contained in a discrete
subring $\Gamma$ of $\C$, and 

ii) there exists $\alpha>0$ such 
that $|z_n|\geq 1+\alpha$ for all $n\geq 1$. 

\noindent Then   
$q_n\neq 0$ for all $n\geq 0$, and $\frac {p_n}{q_n} \to z$ as $n\to \infty$. 
Also, for all $n$ such that
$|q_{n-1}|<|q_n|$,  $|z-\frac{p_n}{q_n}|\leq \alpha^{-1}|q_n|^{-2}.$
\end{proposition}

\proof We note that since $\Gamma$ is a discrete subring of $\C$, 
for any $p\in \Gamma \backslash \{0\}$ we have $|p|\geq 1$. 
Now if $q_n=0$, for some $n\geq 1$, then 
by Proposition~\ref{prop}\,(i) we should have $|p_n|=|z_1\cdots
z_{n+1}|^{-1}\in (0,1)$, which is not possible since $p_n\in \Gamma$. 
Hence $q_n\neq 0$ for all $n\geq 0$. Since $q_n\in \Gamma$ this 
implies that $|q_n|\geq 1$ for all $n$. Therefore, $|z-\frac{p_n}{q_n}|=
|q_n|^{-1}|q_nz-p_n|= |q_n|^{-1}||z_1\cdots z_{n+1}|^{-1}\leq 
(1+\alpha)^{-n}\to 0$, and hence $\frac{p_n}{q_n} \to z$ as $n\to \infty$. 
When $|q_{n-1}|<|q_n|$, by Proposition~\ref{prop} 
$|z-\frac{p_n}{q_n}|\leq |q_n|^{-2}(|z_{n+1}|-|\frac{q_{n-1}}{q_n}|)^{-1}
\leq \alpha^{-1}|q_n|^{-2},$ since $z_n\geq 1 +\alpha$. 
 \hfill $\Box$

\bigskip
A standard way to generate iteration sequences is via algorithms. 
Let $\Lambda$ be a countable subset of $\C$ such that for every $z\in \C$ 
there exists $\lambda \in \Lambda$ such that $|z-\lambda|\leq 1$. 
By a $\Lambda$-valued {\it algorithm} we mean a map $f:\C \to \Lambda$
such that for all $z\in \C$, $|z-f(z)|\leq 1$. 
Let $K$ denote the subfield of $\C$ generated by $\Lambda$; we note 
that $K$ is also countable. For any $z\in \C\backslash
K$ a $\Lambda$-valued  algorithm $f$ as above yields  an 
iteration sequence defined by $z_0=z$ and 
$z_{n+1}=(z_n-f(z_n))^{-1}$ for all $n\geq 0$; for $z\in \C\backslash K$, 
it may be
observed successively that all $z_n\in \C\backslash K$ and hence 
$z_n\neq f(z_n)$, so  $z_n-f(z_n)\neq 0$. 

\begin{definition}
We call  the set 
$\{z-f(z)\mid z\in \C\backslash K\}$ the {\it fundamental set} of the 
algorithm 
$f$. 
\end{definition}

When $\Lambda$ is a discrete subring of $\C$ we have the following.

\begin{theorem}\label{prop3}
Let $\Gamma$ be a discrete subring of $\C$ and let $f:\C\to \Gamma$ be
a $\Gamma$-valued algorithm such that the fundamental set of $f$ is
contained in a ball of radius $r$ centered at $0$, where
$0<r<1$. Let  $K$ be the subfield
generated by $\Gamma$.  Let $z\in \C\backslash K$ and let
$\{z_n\}_{n=0}^\infty$ be 
the iteration sequence for $z$ with respect to $f$.  Let 
$\{a_n\}_{n=0}^\infty$ be the associated sequence of partial
quotients, and  $\{p_n\}, \{q_n\}$ the corresponding $\cal Q$-pair. Then   

i) $q_n\neq 0$ for all $n\geq 0$, and $\frac {p_n}{q_n}\to z$ as $n\to \infty$, 
and

ii) for every $n$ such that 
$|q_{n-1}|<|q_n|$ we have  $|z-\frac{p_n}{q_n}|\leq \frac r{1-r}|q_n|^{-2}.$

 \end{theorem}

\proof Under the condition as in the hypothesis 
$|z_n-a_n|\leq r$, for all $n\geq 0$. Hence for all $n\geq
1$ we have $|z_n|=|z_{n-1}-a_{n-1}|^{-1} \geq r^{-1}$. Thus condition~(ii) of  
Proposition~\ref{prop2} holds, with $\alpha =r^{-1}-1$, and hence the
theorem follows from the proposition. \hfill $\Box$

\medskip
When $\Lambda$ is a discrete subset  we have
an algorithm $f$ arising canonically, where we choose, for $z\in
\C$, $f(z)$ to be the element of $\Lambda$ nearest to $z$;
the map is defined uniquely by this only for $z$ in the  complement a
countable set of lines  
(consisting of points which are equidistant from two distinct points of 
$\Lambda$), but we consider it extended to $\C$ through
some convention - the specific choice of the extension will not play 
any role in our discussion. We call this the {\it nearest element
  algorithm} with respect to $\Lambda$; when $\Lambda$ is a ring 
of ``integers'', such as the Gaussian or Eisenstein integers, the algorithm 
will 
be referred to as the nearest integer algorithm of the corresponding 
ring.

\begin{remark}
{\rm  
It can be seen that any discrete subring  $\Gamma$ of $\C$ (containing~$1$), 
other than $\Z$, has the form $\Z[i\sqrt k]$ or $\Z[\frac 12+\frac i2 
\sqrt{4l-1}]$,
with $k, l\in \N$. From among these,  the requirement 
that there be an element of $\Gamma$ within distance $1$ from 
every $z$ in $\C$ (enabling continued fraction expansions to be defined for all 
$z\in \C$) is 
met for $\Z[i\sqrt k]$, $1\leq k \leq 3$, and 
$\Z[\frac 12+\frac i2 \sqrt{4l-1}]$, $1\leq l \leq 3$;  
for $k=1$ and $l=1$  these are the rings of Gaussian integers and  Eisenstein 
integers respectively.  With respect to the nearest integer algorithm the 
fundamental set is the square with vertices at $\pm \frac 12+\pm 
\frac {\sqrt k} 2 i$ for $\Gamma =\Z[i\sqrt k]$, $k=1,2,3$,  and 
for $\Gamma =\Z[\frac 12 +\frac i2 \sqrt \tau]$, with $\tau =3,7$ or 
${11}$ it  is a hexagon (not regular in  the last
two cases) with 
vertices at $\displaystyle{\pm \frac 12 \pm \frac {\tau-1}{4\sqrt \tau}i}$ and 
$\displaystyle{\pm \frac { \tau+1}{4\sqrt \tau}i}$ respectively; thus the vertices
lie on the circle, centered at the origin, with radius $\frac 12 
\sqrt {(1+k)}$, 
$k=1,2,3$ in the former case, and  $\displaystyle \frac { \tau+1}{4\sqrt 
\tau}$, with 
$\tau =3,7,{11}$, in the latter case; thus the fundamental set is contained in  
a  the open unit ball, except for  $\Z[i\sqrt 3]$. 
Thus, except when $\Gamma=\Z[i\sqrt 3]$  (a case
not considered in literature), by Theorem~\ref{prop3}, $q_n\neq 0$ for all 
$n\geq 0$ and
 ${\displaystyle{\frac {p_n}{q_n}\to z}}$ as $n\to \infty$, for the 
continued fraction 
expansion with respect to the respective nearest integer algorithms. 
}
\end{remark}

\begin{remark}
{\rm 
The second assertion in  
Theorem~\ref{prop3} highlights the usefulness of  establishing 
the monotonicity of  $\{|q_n|\}$,  to complete the 
picture; the monotonicity condition will  also be involved in proving 
the analogue of the Lagrange theorem (see~Corollary~\ref{Lagr-algor}). 
The latter was proved by Hurwitz for the nearest integer
algorithms with respect to the rings of Gaussian integers and Eisenstein 
integers.  It  was proved by Lund for the nearest 
integer algorithm on $\Z[i\sqrt 2]$, as noted in \cite{L}, where it is also 
stated, without proof, that monotonicity holds for $\Z[\frac 12+
\frac i2 \sqrt \tau ]$,
$\tau =3, 7$ or $11$, for the nearest integer algorithm as well as another
 variation
 of it (in each case; see \cite{L} for details). These verifications involve 
elaborate arguments involving ``succession rules''; namely certain restrictions
that hold for the succeding partial quotient in the expansion. 
In \cite{DN-compl} we 
established monotonicity for a variety of algorithms with values in the ring
of Gaussian integers, under a general condition. In the following section we 
extend the idea and  introduce a condition
on the partial quotients which ensures such monotonicity independent of 
the algorithm involved, and even the domain for drawing the partial quotients.  }
\end{remark}

\section{Monotonicity of the denominators of the convergents}

In this section we describe certain general conditions which ensure 
that the denominators of the convergents grow monotonically in size, viz. 
$|q_{n+1}|> |q_n|$ for all $n\geq 0$ in the notation as above.

For $z\in \C$ and $r >0$ we denote by $B(z,r)$ and $\bar B(z,r)$, 
respectively
the open and closed balls with center at $z$ and radius $r$. We note
that if $|z|>r$ then  $\bar B(z,r) \subset \C^*$ and  the sets 
$B(z,r)^{-1}$ and $\bar B(z,r)^{-1}$
(consisting of the inverses of elements from the respective sets) are 
given by $\displaystyle{{B\left (\frac{\bar z}
{|z|^2-r^2},\frac{r}{|z|^2-r^2}\right )}} $ and 
$\displaystyle{{\bar B\left (\frac{\bar z}
{|z|^2-r^2},\frac{r}{|z|^2-r^2}\right )}} $ respectively.   

\smallskip
\begin{definition}
{\rm A sequence $\{a_n\}_{n=0}^\infty$ in $\C$ is said to satisfy 
{\it Condition $\cal C$}
if $|a_n|>1$ for all $n\geq 1$ and whenever $|a_{n+1}|<2$ for some $n\geq 1$ 
then we have
 $|(|a_{n+1}|^2-1)a_{n}+\bar a_{n+1}| \geq |a_{n+1}|^2 $.
}
\end{definition}
\begin{theorem}\label{prop:mono}
Let $\{a_n\}_{n=0}^\infty$ be a sequence in $\C$ satisfying Condition~$\cal C$
and let 
$\{p_n\}, \{q_n\}$ be the corresponding $\cal Q$-pair.
Then  $|q_{n+1}|>|q_n|$ for all $n\geq 1$.  
\end{theorem}

\proof Suppose, if possible, that there exists $n\geq 1$ such 
that $|q_{n+1}|\leq |q_n|$, and let $m\geq 1$ be the smallest such number. 
Thus we have $|q_{m+1}|\leq |q_{m}|$ and  $|q_{n+1}|> |q_{n}|$ for 
$n=1,\dots, m-1$. In particular $q_n\neq 0$ for $n=1, \dots, m$. 
For all $0\leq n\leq m$ let $r_n=q_{n+1}/q_n$; then $r_n>1$ for $n=0,1,\dots, 
m-1$, and $r_m\leq 1$.
From the recurrence relations for $\{q_n\}$ we have 
$r_n=a_{n+1}+r_{n-1}^{-1}$ for all $1\leq n\leq m$. 
In particular $r_{m-1}^{-1}\in 
\bar B(-a_{m+1},|r_m|)\subset  \bar B(-a_{m+1},1) $ and, since $|a_{m+1}|>1$, 
this implies $r_{m-1}\in 
\displaystyle{\bar B\left (\frac{-\bar a_{m+1}}{| a_{m+1}|^2-1},
\frac 1{| a_{m+1}|^2-1}\right )}$. We have $r_{m-1}=a_{m}+r_{m-2}^{-1}$, 
and together with the preceding  
conclusion we get that $a_m\in \displaystyle{\bar 
B\left (\frac{-\bar a_{m+1}}{| a_{m+1}|^2-1}, |r_{m-2}^{-1}|+
\frac 1{| a_{m+1}|^2-1}\right )}$.
 In turn, since  $|r_{m-2}|>1$, $a_m$ is contained in the open ball
$\displaystyle{B\left (\frac{-\bar a_{m+1}}{| a_{m+1}|^2-1}, 1+\frac 1
{| a_{m+1}|^2-1}\right )} $. Thus $$|(| a_{m+1}|^2-1)a_m +\bar a_{m+1}|<
(| a_{m+1}|^2-1)+1 =|a_{m+1}|^2.$$ 
On the other hand, since 
 $r_m=a_{m+1}+r_{m-1}^{-1}$ we have $|a_{m+1}|\leq |r_m|+|r_{m-1}^{-1}|<2$.
Together with the above conclusion this contradicts the 
condition in the hypothesis. Therefore $r_n>1$ for all $n\geq 0$,
or equivalently  $|q_{n+1}|>|q_n|$ for all $n\geq 0$. This proves the 
proposition.  \hfill $\Box$

\begin{remark}
{\rm Let $\Gamma$ be a discrete subring of $\C$ and $f:\C\to \Gamma$ be 
a $\Gamma$-valued algorithm such that the fundamental set of $f$ is 
contained in a ball of radius $0<r<1$. 
Let $z\in \C^*\backslash K$, where $K$ is the subfield generated by 
$\Gamma$, and let  $\{a_n\}_{n=0}^\infty$ be the sequence of partial 
quotients for $z$ with respect to $f$, and  $\{p_n\}, \{q_n\}$ 
be the   $\cal Q$-pair corresponding to  $\{a_n\}_{n=0}^\infty$.
If  $\{a_n\}_{n=0}^\infty$ satisfies Condition~$\cal C$, 
then by Theorem~\ref{prop:mono} $|q_{n+1|}>|q_n|$ for all 
$n\geq 0$, and  by Theorem~\ref{prop3} $\displaystyle |z-\frac 
{p_n}{q_n}|\leq c|q_n|^{-2}$ for all $n\geq 0$, 
with $c=\frac r{1-r}$. From a Diophantine point of view these are
only weak estimates - but seem to be of significance on account of  
generality of their context. 
In \cite{L} optimal values for such a constant $c$ 
are described for continued fraction expansions with respect to the nearest 
integer algorithms, and also a variation in the case of   
$\Z[\frac 12+\frac i2 \sqrt{\tau }]$, $\tau =3,7$ or $ 11$. It would be 
interesting to know 
similar optimal values for more general algorithms. 
 }
 \end{remark}

\begin{remark}{\rm 
{Let $\frak G$ denote the ring of 
Gaussian integers, viz. $\frak G=\Z[i]$. Let $z\in \C$ and  $\{z_n\}_{n=0}^\infty$ 
be an iteration sequence for $z$ such that 
$a_n=z_n-z_{n+1}^{-1} \in \frak G$ for all $n\geq 0$. 
For $a\in \frak G$, $1<|a|<2$ if and only if $a=\pm 1\pm i$, or equivalently  
$|a|=\sqrt 2$. Thus in this case  Condition~$\cal C$  
reduces to that for all $n\geq 1$, $|a_n|>1$ and either $|a_{n+1}|\geq 2$ or 
$|a_n+ \bar a_{n+1}|\geq 2$. This corresponds to Condition~(H') in  
 \cite{DN-compl}, used for obtaining a conclusion as in Theorem~\ref{prop:mono}
as above; a special case of  Theorem~\ref{prop:mono} was obtained there, 
in Theorem~6.11, only after proving other results about the asymptotic
growth of $|q_n|$'s. It may also be recalled here that the
sequence $\{a_n\}$ obtained by application of the nearest (Gaussian) integer 
algorithm, starting with a $z\in \C\backslash \Q(i)$ 
may not satisfy Condition~$\cal C$ (the second part) (see 
\cite{DN-compl}, \S~5 for details). The sequences 
corresponding to the nearest 
integer algorithm satisfy a weaker condition, named Condition~(H) in 
\cite{DN-compl}, which also suffices to obtain the conclusion as in 
Theorem~\ref{prop:mono}; the condition however is rather technical and not
amenable to generalization. In \cite{DN-compl} another algorithm, named
PPOI (acronym for partially preferring odd integers) was 
introduced, producing a continued fraction expansion in terms of Gaussian 
integers for which Condition~(H') is satisfied. We shall however show 
in the following sections that in the case of the Eisenstein integers 
the sequences
corresponding to the nearest integer algorithm, as also certain other 
algorithms, satisfy Condition $\cal C$.}}

\end{remark}

\section{Lagrange theorem for continued fractions}

In this section we prove an analogue of the classical Lagrange theorem, 
about the continued fraction expansion
being eventually periodic if and only if the number is a quadratic
surd.  We shall continue to follow the notation as before. 

Let $K$ be a subfield of $\C$. 
A number $z\in \C$ is called a {\it quadratic surd} over $K$ if $z\notin K$ 
and it is a root of a quadratic polynomial over $K$. 

\begin{proposition}\label{dir}
Let  $z\in \C\backslash K$ and
$\{z_n\}_{n=0}^\infty$ be an iteration  sequence for $z$ such that 
such that $|z_n|>1$ for all $n\geq 1$. Let 
$a_n=z_n-z_{n+1}^{-1}$,  $n\geq 0$, be the corresponding sequence of 
partial quotients and suppose that $a_n$, $n\geq 0$, are all contained 
in a discrete subring $\Gamma$ 
of $\C$ contained in $K$. Let  $\{p_n\}$,
$\{q_n\}$ be the   corresponding $\cal Q$-pair.  
If $z_m=z_n$ for some  $0\leq m< n$, then $z$ is a quadratic
surd over $K$. 
\end{proposition}

\proof Clearly, for all $m\geq 0$, $\{z_{m+k}\}_{k=0}^\infty$ is an iteration 
sequence for $z_m$
and $z$ is a quadratic surd if and only if $z_m$ is a 
quadratic surd.  Hence in proving the proposition we may assume that 
$z_m=z$, or equivalently that 
$m=0$. Let $n\geq 1$ be such that $z_n=z$. 
By Proposition~\ref{prop} we have $(q_{n-1}z-p_{n-1})z_n=(-1)^{n-1}(z_1
\cdots z_n)^{-1}z_n=(q_{n-2}z-p_{n-2})$. Since by hypothesis $z_n=z$, 
we get $q_{n-1}z^2-(p_{n-1}+q_{n-2})z+p_{n-2}=0$. 
Suppose, if possible, that $q_{n-1}=0$. 
Then   $|p_{n-1}|=|q_{n-1}z-p_{n-1}|=|z_1 \cdots z_n|^{-1} \in (0,1)$, 
which is not possible since $p_{n-1}$ is contained in a
discrete subring $\Gamma$ of $\C$.  Thus $q_{n-1}\neq 0$, and  we 
see that $z$ satisfies a quadratic polynomial over $K$. 
Since $z\notin K$ it follows that  $z$ is a quadratic surd over $K$. 
 \hfill $\Box$

\medskip
We now prove the following converse of this. The proof follows what is 
now a standard strategy (cf. \cite{HW} for instance) for proving such a 
result, with variations in 
the hypothesis; the main purpose here is to bring out  
a general formulation which at the same time is focused enough 
and amenable to a
brief treatment.

\begin{theorem}\label{thm:lagr}
Let $\Gamma$ be a discrete subring of $\C$ and $K$ be the quotient field of 
$\Gamma$. Let $z$ be a quadratic surd over $K$. 
Let $\{z_n\}_{n=0}^\infty$ be  an  
iteration sequence for $z$ such that  the corresponding 
sequence  $\{a_n\}_{n=0}^\infty$ of partial quotients is contained in $\Gamma$. 
Let   $\{p_n\}$, $\{q_n\}$ be the $\cal Q$-pair corresponding to  
$\{a_n\}_{n=0}^\infty$.
Suppose that the following conditions are satisfied:

i) there exists $\alpha >0$ such that $|z_n|>1+\alpha$ for all $n\geq 1$; and

ii)  $|q_{n}|\to \infty $ as   $n\to \infty$.  

\noindent Then the set $\{z\in \C\mid z=z_n \mbox { for some } n\}$ is 
finite.  Consequently, if  $\{z_n\}_{n=0}^\infty$ is an iteration sequence
associated with an algorithm then $\{a_n\}_{n=0}^\infty$ is eventually
periodic.  

\end{theorem}

\noindent{\it Proof}:
Let   $a,b,c \in K$, with $a \neq 0$, be such that $ az^2+bz+c=0$. Since 
$K$ is the quotient field of $\Gamma$ we may without loss of generality assume 
that $a,b,c\in \Gamma$. 
By Proposition~\ref{prop}(iii) we have $z=\displaystyle{\frac{z_{n+1}p_{n}+p_{n-1}}
{z_{n+1}q_{n}+q_{n-1}}}$, for all $n\geq 0$, and hence 
$$
a\left (\frac{z_{n+1}p_{n}+p_{n-1}}{z_{n+1}q_{n}+q_{n-1}}\right )^2+
b\left (\frac{z_{n+1}p_{n}+p_{n-1}}{z_{n+1}q_{n}+q_{n-1}}\right )+c=0.
$$
For all $n\geq 0$ let  
$$
A_n=ap_{n}^2+bp_{n}q_{n}+cq_{n}^2,\, C_n=A_{n-1},  \hbox{\rm and}$$  
$$
B_n=2ap_{n}p_{n-1}+b(p_{n}q_{n-1}+q_{n}p_{n-1})+2cq_{n}q_{n-1}.  
$$ 
Then $A_n, B_n,C_n\in \Gamma$, for all $n$, and the above equation can be 
readily simplified to $A_nz_{n+1}^2+B_nz_{n+1}+C_n=0.$
The polynomial $a\zeta ^2+b\zeta+c$ has no root in $K$, and hence it 
now follows that $A_n\neq 0$ for all $n$. Now, we have
$A_n = (ap_{n}^2+bp_{n}q_{n}+cq_{n}^2)-q_{n}^2(az^2+bz+c)$, and 
the latter expression can be rewritten as 
$(p_{n}-zq_{n})(a(p_{n}-zq_{n})+(2az+b)q_{n}).$ Under the
conditions in the hypothesis, by Proposition~\ref{prop2} we have 
$|q_nz -p_n|\leq 
\alpha^{-1} |q_n|^{-1}$. Therefore by substitution we get that, we get that
$$|A_n|\leq  \alpha^{-1} |q_n|^{-1}(a \alpha^{-1} |q_n|^{-1}+|2az+b||q_{n}|)
= \alpha^{-1}|2az+b|+\alpha^{-2} a|q_n|^{-2}.$$
Since by hypothesis $|q_n|\to \infty$ as $n\to \infty$, the above 
observation implies that  $\{A_n\mid n\geq 0\}$ is a
bounded set, and since  $A_n\in \Gamma$ for all $n$ it further follows 
that $\{A_n\mid n\geq 0\}$ is finite. 
 Since $C_n=A_{n-1}$ for all $n\geq 1$, we also have  $\{C_n\mid n\geq 0\}$
is finite. An easy  computation shows  that for all $n\geq 0$,
$B_n^2-4A_nC_n= b^2-4ac$. It follows that $\{A_n\zeta ^2+B_n\zeta +C_n \mid 
n\geq 0\}$ is a
finite collection of  polynomials. Since each $z_n$ is a root of one 
of these polynomials we get that  $\{z_n\mid n\geq 0\}$  is finite. This
proves the first statement in the theorem. 

Now suppose that  $\{z_n\}_{n=0}^\infty$ is an iteration sequence
associated with an algorithm. By the first part we get that there
exist $m\geq 0$ and $k\geq 1$ such that $z_{m+k}=z_m$. Since 
$\{z_n\}_{n=0}^\infty$ are determined algorithmically, this implies
that $z_{n+k}=z_n$ for all $n\geq m$. In turn we get that  $a_{n+k}=a_n$ 
for all $n\geq m$, that is, $\{a_n\}_{n=0}^\infty$ is eventually periodic. 
This proves the theorem. \hfill $\Box$ 

\medskip
The following Corollary (together with Proposition~\ref{dir}) gives a 
generalisation of the classical Lagrange theorem of quadratic irrationals. 

\begin{corollary}\label{Lagr-algor}
Let $\Gamma$ be a discrete subring of $\C$ and $K$ be the quotient field of 
$\Gamma$. Let $f:\C\to \Gamma$ be a $\Gamma$-valued algorithm such that the 
fundamental set of $f$ is contained in a ball of radius $0<r<1$. 
Let $z$ be a quadratic surd over $K$. 
Let  $\{a_n\}_{n=0}^\infty \subset \Gamma$ be the sequence of partial quotients 
with respect to $f$ and let   $\{p_n\}$, $\{q_n\}$ be the corresponding 
$\cal Q$-pair. Suppose that $\{a_n\}$ satisfies Condition~$\cal C$. 
Then $\{a_n\}_{n=0}^\infty$ is eventually
periodic.

\end{corollary}

\proof This follows from Theorem~\ref{thm:lagr}:  Condition (i) in the 
hypothesis
of the theorem is satisfied since the fundamental set of $f$ is contained a 
ball of radius $r<1$. Under Condition~$\cal C$ as in the hypothesis by 
Theorem~\ref{prop:mono} $\{|q_n|\}$ is strictly monotonically increasing and, 
since $\Gamma$ is a discrete subring, we get that $|q_n|\to \infty$ 
as $n\to \infty$,
so condition~(ii) holds. Hence the Corollary. 
\qed

\section{Continued fractions for Eisenstein integers}

We shall now apply the results of the preceding sections to a class of 
algorithms with values in the ring of Eisenstein integers. Let $\frak E$ 
be the ring of  Eisenstein integers in $\C$, 
viz. $\frak E=\{x+y\omega \mid x, y\in \Z\}$, where $\omega$ is a primitive 
cube root 
of unity, which we shall realise as $-\frac 12+\frac {\sqrt 3}2 i$. 
Let 
$\rho =\frac 12+\frac {\sqrt 3}2i$ (which is a primitive $6$th root of unity). 
Then $\rho =\omega +1$, and every $z\in \frak E$ can also be expressed 
as $x+y\rho$, with  $x,y 
\in \Z$. For convenience we shall also use the notation $j$ for $\sqrt 3\,i$. 
Then every $z\in \frak E$ can be expressed as $\frac 12(x+yj)$ with 
$x,y\in 2\Z$, viz. $x+y$ an even integer. We shall write the~$6$~th roots of 
unity as $\rho^k$, with $k\in \Z$, the integer $k$ being understood to be 
modulo~$6$. 

Given a $\frak E$-valued algorithm $f$ we shall 
denote by $\Phi_f$ its fundamental 
set of $f$, and by $C_f(a)$, for $a\in \frak E$, the set 
$\{z\in \C\mid f(z)=a\}$. 

\begin{theorem}\label{thm:eisen}
Let $\frak E$ be the ring of Eisenstein integers and let  
$f:\C\to \frak E$ be a $\frak E$-valued algorithm and let $\Phi=\Phi_f$.  
Suppose that 

a) $\Phi \subset B(0,r)$, for some $0<r<1$, 

b) $|f(\zeta)|>1$ for all $\zeta \in \Phi^{-1}$, and 

c) for $0\leq k\leq 5$ and $t\in \{-1+j, j,1+j\}$, the sets
 $ \rho^{-k} t+(C_f(\rho^kj))^{-1}$ and  $C_f( \rho^{-k}t)\cap \Phi^{-1}$ are disjoint. 

\noindent Let $K$ be the subfield generated by $\frak E$. 
Let $z\in \C\backslash K$, $\{a_n\}_{n=0}^\infty$ be the sequence of partial 
quotients of $z$ corresponding to the algorithm $f$ and  $\{p_n\}$, 
$\{q_n\}$ be 
the $\cal Q$-pair corresponding to $\{a_n\}_{n=0}^\infty$. Then the following 
conditions are satisfied:

i) $|q_n|>|q_{n-1}|$ for all $n\geq 1$, and in particular $q_n\neq 0$ 
for all $n$. 

ii)  $\frac {p_n}{q_n}\to z$ as $n\to \infty$ and moreover $|z-\frac {p_n}{q_n}|
\leq \frac r{1-r}|q_n|^{-2}$ for all $n$. 

iii) $z$ is a quadratic surd over $K$ if and only if $\{a_n\}$ is eventually
periodic.

\noindent
\end{theorem}

\proof Let the notation be as in the hypothesis. Also let 
 $\{z_n\}_{n=0}^\infty$ denote the 
corresponding iteration sequence for $z$ with respect to $f$.
Since by Condition~(b) $|f(\zeta)|>1$ for all $\zeta \in \Phi^{-1}$ 
it follows that $|a_n|>1$
for all $n\geq 1$. We shall show that $\{a_n\}$ satisfies Condition~$\cal C$. 
For this we first note that for $a\in \frak E$, if $1<|a|<2$ then 
$|a|=\sqrt 3$, and $a=\rho^kj$ for some $k\in \Z$. Hence we need to show 
that for $n\geq 1$ if $a_{n+1}=\rho^kj$, $k\in \Z$, 
then $|2a_n-\rho^{-k}j|\geq 3$. Let if possible  $n\geq 1$ be such that 
$a_{n+1}=\rho^kj$, $k\in \Z$, and $|2a_n-\rho^{-k}j|< 3$. We write 
$a_n$ as $\frac 12\rho^{-k}(x+yj)$, with $x,y\in 2\Z$. Then by 
the above condition we have $3>|2a_n-\rho^{-k}j|= |2a_n\rho^{k}-j|=|x+yj-j|$, 
and hence $x^2+3(y-1)^2<9$. Also, since $|a_n|>1$ we have $x^2+3y^2\geq 12$. 
The only common solutions to this, with $x+y$ even, are $x=0$ or 
$\pm 2$  with $y=2$. Thus $a_n\in \rho^{-k} \{-1+j, j,1+j\}$.  
We have $z_n\in \Phi^{-1}$ (as $n\geq 1$), $z_n\in C_f(a_n)$, and also 
$z_n=a_n +z_{n+1}^{-1}\in a_n +(C_f(\rho^kj))^{-1}$. Since 
$a_n\in \rho^{-k}\{-1+j, j,1+j\}$  this contradicts  
condition~(b) in the hypothesis, for $t=\rho^ka_n$. Hence the desired condition 
as above holds for all $n\geq 1$. Assertion~(i) in the theorem now 
follows from  Theorem~\ref{prop:mono},
and, together with Condition~(a) in the hypothesis, it implies 
assertions (ii) and~(iii),  in view of 
 Theorem~\ref{prop3} and Corollary~\ref{Lagr-algor} respectively. \qed

\begin{corollary}\label{cor:eisen}
Let $\frak E$ be the ring of Eisenstein integers and let  
$f:\C\to \frak E$ be a $\frak E$-valued algorithm such that the following 
conditions hold:

a) $C_f(a)$ is contained in $B(a,\frac12 (\sqrt 5-1))$ for all 
$a\in \frak E$, and 

b) for $0\leq k\leq 5$, $C_f(\rho^kj)$ is contained in  
$B(\rho^kj,\sqrt \lambda) $, where $\lambda = \frac 14(5-\sqrt {13})$. 

\noindent Then  statements (i), (ii) and (iii) as in Theorem~\ref{thm:eisen} 
are satisfied.  In particular they are satisfied for the nearest integer 
algorithm. 
\end{corollary}

\proof Condition~(a) as in Theorem~\ref{thm:eisen} is evidently satisfied 
for any $f$ as above. 
 We show that conditions~(b) and~(c) are also satisfied. By condition~(a)
in the hypothesis $\Phi_f \subset B(0,r)$, for $r=\frac12 (\sqrt 5-1)$. 
Hence for $\zeta \in \Phi_f^{-1}$ we have $|\zeta|
>r^{-1}= 1+r$, and 
since $f(\zeta) \in B(\zeta , r)$, this shows that $|f(\zeta)|>1$,
 thus proving condition~(b). 

Now let $0\leq k\leq 5$ and $t\in \{-1+j, j,1+j\} $ be given. To begin with 
consider any $r>0$  such that $C_f(\rho^kj) \subset B(0,r)$; we shall show 
that  condition~(c) of Theorem~\ref{thm:eisen} holds when $r< \sqrt \lambda$.  
Putting $\sigma =(3-r^2)^{-1}$ (as temporary notation for 
convenience), we have 
$$ \rho^{-k}t+
(C_f(\rho^kj))^{-1} \subset  \rho^{-k}t+B(\rho^{k}j, r)^{-1}= \rho^{-k}t+ \rho^{-k} 
B(-\sigma j, \sigma r),$$
which is the same as
$\rho^{-k} B(t-\sigma j, \sigma r).$ 
 Since $\Phi_f^{-1}$ is complimentary to $B(0,r^{-1})$, to 
prove condition~(c) it now suffices to show that,  for all 
$t\in \{j-1,j,j+1\}$, 
$B(t-\sigma j, \sigma r)$ is contained in $B(0,r^{-1})$; the condition is now 
independent of $k$. For $t=j$ it suffices to note that $|t-\sigma j|+ \sigma r=
(1-\sigma)\sqrt 3 +\sigma r=\sqrt 3 -\sigma (\sqrt 3 -r)=\sqrt 3 
-(\sqrt 3 +r)^{-1}$, on substituting for $\sigma$. 
The last expression is less than $r^{-1}$ when $\sqrt 3 r^2+r-\sqrt 3<0$, 
viz. if  $r< (\sqrt {13} -1)/2\sqrt 3 \approx 0.752\cdots $, so it holds 
in particular 
for $r<\sqrt {\lambda } \approx 0.590\cdots $ as in the hypothesis.  

It remains to consider the case of $t=\pm 1+j$, and by symmetry it suffices 
to consider the case $t=1+j$. We need to verify that
$|1+(1-\sigma)j|+\sigma r<r^{-1}$, or equivalently $1+3(1-\sigma)^2<
(r^{-1} -r\sigma)^2=r^{-2} 
(1-r^2\sigma)^2$.  Substituting for  $\sigma$ as $(3-r^2)^{-1}$ and eliminating
the denominators the condition reduces to 
$$(3-r^2)^2r^2+3(2-r^2)^2r^2-(3-2r^2)^2<0.$$ 
Let $s=r^2$ and $P$ be the polynomial $P(s)=4s^3-22s^2+33s-9$; the 
above expression then coincides with $P(r^2)$. Now we see that $$P(s)=(s-3)
(4s^2-10s+3)=(s-3)(s-\lambda)(s-\mu),$$ where $\lambda$ is as in the 
hypothesis and $\mu= \frac 14 (5+\sqrt {13})$ is its quadratic conjugate.
Hence $P(s)<0$ for $s<\lambda$, and hence $P(r)<0$
for $r<\sqrt \lambda$. Hence Condition~(c) holds and therefore, by 
Theorem~\ref{thm:eisen} the assertions (i), (ii) and (iii) 
as in the theorem hold.  For the nearest integer algorithm the fundamental
set is a regular hexagon contained in $\bar B(0,1/\sqrt 3) \subset B(0, \sqrt 
\lambda)$ and so the assertions hold, as a particular case. \qed
\medskip

\begin{example}
{\rm
Let   $P$ denote the closed parallelogram  with 
vertices at  $0,1,\rho$ and $1+\rho$. Then $\C$ 
is tiled by $\{a+P\}_{a\in \frak E}$ and it suffices to define the algorithm on 
each $a+P$, $a\in P$; the points on the boundaries may be assigned a specific 
tile $a+P$ by some convention in applying the following. Let $0<r<1$ and 
$V= \{0,1,\rho, 1+\rho\}$. 
Let $P_v$, $v\in V$, be disjoint subsets of $P$ such that $P_v\subset B(v,r)$ 
and $P=\cup_{v\in V} P_v$. It may be seen that such partitions exist for 
$r>1/\sqrt 3$. 
Then we can define an algorithm $f:\C\to \frak E$, 
by setting, for any $a\in \frak E$ and $\zeta \in P_v$, $f(a+\zeta)= a+v$. 
(We note that the choice of the partition as above may also be made dependent 
on $a$).   Then $C_f(a)\subset B(a,r)$ for all $a\in \frak E$. If we choose 
$r\leq \frac 12 \sqrt {(5-\sqrt {13})}$, then the conditions in 
Corollary~\ref{cor:eisen} are satisfied, and therefore the statements as
in the conclusion hold for such an algorithm. } 
\end{example}

\section{Exponential growth of $\{|q_n|\}$}

It is known in the case of various algorithms over the ring of Gaussian 
integers that the sequence $\{|q_n|\}$ increases exponentially; see 
\cite{DN-compl}. We shall show that analogous assertion also holds in the 
case of Eisenstein integers. For simplicity we shall restrict to the nearest 
integer algorithm in this respect; extension to some of the algorithms as
 the second half of Theorem~\ref{thm:eisen}, seems feasible but involves 
some cumbersome computations, which do not seem worthwhile for the present. 

\begin{theorem}\label{thm:mono} 
Let $\frak E$ be the ring of Eisenstein integers, $K$ the subfield generated by 
$\frak E$. Let $z\notin K$ and let $\{a_n\}$ be the sequence of partial
 quotients of $z$ corresponding to the nearest integer algorithm. Let 
 $\{p_n\}, \{q_n\}$ be  $\cal Q$-pair corresponding
to $\{a_n\}$. Then $\displaystyle{|\frac{q_{n+1}}{q_{n-1}}|>\frac 32}$ for all $n\geq 1$. 
\end{theorem}

For this we first prove the following. 

\begin{proposition}\label{prop5} Let the notation be as in 
Theorem~\ref{thm:mono}. Then for all  $n\geq 1$ we have the following:

i)  if $a_n=j\rho^k$, $k\in \Z$, then 
$a_{n+1}\rho^k=\frac 12 (x+yj)$ with $x,y\in 2\Z$,  such 
that  $|\frac 12 x|\leq 2-\frac 32 y$.

ii)  if $a_n=2\rho^k$, $k\in \Z$, then $a_{n+1}\rho^k=
\frac 12 (x+yj)$, with $x,y\in 2\Z$, and $x\geq -2$.

\end{proposition}

\proof Let $\{z_n\}$ be the iteration sequence of $z$ (with respect to the 
nearest integer algorithm). Let $H$ be the hexagon with vertices at 
$\frac 13 \rho^kj$,  $0\leq k\leq 5\}$, the fundamental set of the 
algorithm. 

Let $n\geq 1$ be such that $a_n=j\rho^k$ for some $k$. 
We have $z_n\in a_n+H$, and since $n\geq 1$
we also have $z_n\in H^{-1}$, and so $z_n\notin \cup_{m\in \Z} B(\rho^m,1)$. 
Hence $z_n-a_n \notin  \cup_{m\in \Z} B(\rho^m-j\rho^k,1)$. 
We have  $ B(\rho^m-j\rho^k,1)= \rho^k B(\rho^{m-k}-j,1)$, and
when $m-k=1$ and $2$, we see that $\rho^{m-k}-j=\rho^{-1}$ and $\rho^{-2}$, 
respectively. Thus we get in particular that  
$(z_n-a_n)\rho^{-k} \notin B(\rho^{-1},1)\cup B(\rho^{-2},1) $. Hence 
$z_{n+1}\rho^k=((z_n-a_n)\rho^{-k})^{-1} \notin B(\rho^{-1},1)^{-1}\cup 
B(\rho^{-2},1)^{-1}$. The complements of  $B(\rho^{-1},1)^{-1}$ 
and $B(\rho^{-2},1)^{-1}$ may be seen to be  $\{\sigma +\tau i\mid 
\sigma +\sqrt 3 \tau \leq 1\}$ and  $\{\sigma +\tau i\mid \sigma , 
\tau \in \R, -\sigma +\sqrt 3 \tau \leq 1\}$ respectively ($\sigma$ and $\tau$
understood to be real). 
When $z_{n+1}\rho^k$ belongs to the wedge shaped set consisting of the 
intersection of these two sets, $a_{n+1}\rho^k$ has to belong to the 
intersection of $\{\sigma +\tau i\in \C\mid 
\sigma +\sqrt 3 \tau \leq 2\}$ and  $\{\sigma +\tau i\in \C\mid 
-\sigma +\sqrt 3 \tau \leq 2\}$. With $a_{n+1}\rho^k$ written as 
$\frac 12 (x+yj)$, $x+y $ even, this condition yields $|\frac 12 x|
\leq 2-\frac 32 y$. This proves (i). 

Let $n\geq 1$ be such that $a_n=2\rho^k$ for some $k$. Arguing as above 
we deduce that 
$z_n-a_n\notin \rho^kB(\rho^{m-k}-2,1)$,  for any $m$, and in particular 
choosing $m=k$ we get that $(z_n-a_n)\rho^{-k} \notin B(-1,1)$. Hence 
$z_{n+1}\rho^k\notin B(-1,1)^{-1}$. The complement of $B(-1,1)^{-1}$ is 
the set $\{\sigma +\tau i\mid \sigma \geq -\frac 12\}$, and we see 
that when $z_{n+1}\rho^k$ belongs to it, $a_{n+1}\rho^k$ belongs to 
$\{\sigma +\tau i\mid \sigma \geq -1\}$. Writing $a_{n+1}\rho^k$ as 
$\frac 12 (x+yj)$, $x+y$ even, we get that $x\geq -2$. 
This proves (ii). \qed 

\medskip
In the proof of Theorem~\ref{thm:mono} we use the following simple 
observation, which may be of independent interest. 

\begin{remark}\label{rem2}
{\rm Let $\{a_n\}_{n=0}^\infty$ be a sequence in $\C$ and let 
$\{p_n\}, \{q_n\}$ be the corresponding $\cal Q$-pair. Then  for all  
$n\geq 1$ we have $q_{n+1}=a_{n+1}q_n+q_{n-1}=a_na_{n+1}q_{n-1}+a_{n+1}q_{n-2}
+q_{n-1}$, and hence if  $|q_{n-2}|\leq |q_{n-1}|$ it follows that  
 $$|\frac{q_{n+1}}{q_{n-1}}|=|a_na_{n+1}+1+a_{n+1}\frac{q_{n-2}}
{q_{n-1}}|\geq |a_na_{n+1}+1|-|a_{n+1}|.$$   
}
\end{remark}

\medskip
\noindent {\it Proof of Theorem~\ref{thm:mono}}: 
In view of Remark~\ref{rem2} it would suffice to show that
$|a_na_{n+1} +1|> |a_{n+1}| +\frac 32$, for all $n\geq 1$.
We have $$|a_na_{n+1} +1|- |a_{n+1}| \geq |a_na_{n+1}|-1-|a_{n+1}|=
 (|a_n|-1)|a_{n+1}|-1. $$ 
If $|a_n|>2$ then $|a_n|>\sqrt 7$, and since $|a_{n+1}|\geq \sqrt 3$, we get
$|a_na_{n+1} +1|- |a_{n+1}|\geq (\sqrt 7 -1)\sqrt 3 -1>\frac 32$. 
It remains to consider the cases $|a_n|=\sqrt 3 $ or~$2$.

Suppose  that $|a_n|=\sqrt 3$, so 
$a_n=j\rho^k$, with $k\in \Z$. Then, by Proposition~\ref{prop5} 
we have $a_{n+1}\rho^k=\frac 12 (x+yj)$ with $x,y\in 2\Z$, such 
that  $|\frac 12 x|\leq 2-\frac 32 y$. The last part implies that 
$y\leq 1$, and when $y=1$ 
it further implies, together with $x+y$ being even, that  
$x=\pm 1$, which however is not possible since $|a_{n+1}|\geq \sqrt 3$.  
Hence $y\leq  0$. 
Now, $|a_na_{n+1}+1|=|j\rho^k \cdot \frac 12 (x+yj)
\rho^{-k}+1|=|\frac 12 (xj-3y)+1|$. Therefore  using that $y\leq 0$ we have 
$$|a_na_{n+1}+1|^2=\frac 14 \{3x^2
+(2-3y)^2\}\geq \frac 34(x^2+3y^2)+4=3|a_{n+1}|^2+4. $$ 
We note that  $3|a_{n+1}|^2+4
\geq (|a_{n+1}|+\sqrt {8/3})^2$, as may be seen by considering 
the discriminant of the quadratic difference expression. Thus we have  
$$|a_na_{n+1}+1|\geq |a_{n+1}| +\sqrt {8/3} >|a_{n+1}|+\frac 32,$$ which 
settles the case at hand.  

Now suppose that $|a_n|=2$, so $a_n=2\rho^k$ for some $k\in \Z$. Then by
  Proposition~\ref{prop5}
we have $a_{n+1}\rho^k=\frac 12 (x+yj)$, with $x,y,\in 2\Z$, and 
$x\geq -2$. Hence   $|a_na_{n+1}+1|=|x+yj+1|$. Suppose first that $x\geq 0$. 
 Then $|a_na_{n+1}+1|^2=|x+yj+1|^2=(x+1)^2+3y^2> 4|\frac 12 (x+yj)|^2= 
4|a_{n+1}|^2$. Hence  $|a_na_{n+1}+1|-|a_{n+1}| \geq 2|a_{n+1}|\geq 2\sqrt 3 
> \frac 32$, as desired. 

The only possibilities that remain are $x=-2$ or $-1$. 
We note that since $x^2+3y^2\geq 12$, if $x=-2$ then $|y|\geq 2$ and if $x=-1$
then $|y|\geq 3$. 
Now,  $|a_na_{n+1}+1|^2=|x+yj+1|^2=(x+1)^2+3y^2< 4 y^2$, and   $|a_{n+1}|^2=
|\frac 12(x+yj)|^2=\frac 14 x^2+\frac 34 y^2\leq y^2$. Hence 
$|a_na_{n+1}+1|+|a_{n+1}| <3|y|$.

Suppose $x=-2$. Then  $|a_na_{n+1}+1|^2-|a_{n+1}|^2=1+3y^2-
\frac 14 (4+3y^2)=\frac 94 y^2$, and dividing by the expression estimated 
above, we get  $|a_na_{n+1}+1|-|a_{n+1}|> \frac 3{4}
 |y|\geq \frac 32$, since $|y|\geq 2$, as desired.  Finally 
suppose $x=-1$. Then 
$|a_na_{n+1}+1|^2-|a_{n+1}|^2= 3y^2-\frac 14 (1+3y^2)>2y^2$,
and hence $|a_na_{n+1}+1|-|a_{n+1}|>\frac 23 |y|\geq 2$, since $|y|\geq 3$ in 
this case. Thus we have $|a_na_{n+1}+1|-|a_{n+1}|>\frac 32$ in this case also. 
 This proves the theorem. \hfill $\Box$

\medskip
\begin{remark}\label{rem:const}
{\rm The constant $\frac 32$ involved in
Theorem~\ref{thm:mono} is not  optimal; it was involved closely only in one 
of the special cases in the above argument, where also it can be improved upon
with some detailed 
computations. We shall however not concern 
ourselves here with the aspect of improving it. 
}
\end{remark}

\subsection*{Acknowledgements}
Thanks are due to Lovy Singhal and Ojas
Sahasrabudhe for helpful comments on
an earlier version of this paper. The author is thankful to an anonymous 
referee for pointing out the work of R.B. Lakein~\cite{L} on continued 
fractions 
for complex numbers.

\vskip5mm

\noindent Department of Mathematics\\
Indian Institute of Technology Bombay\\
Powai, Mumbai 400076\\ India

\smallskip
\noindent E-mail: {\tt sdani@math.iitb.ac.in}
\end{document}